\documentclass[11pt,reqno]{amsart}

\usepackage{amssymb}

\catcode`\@=11

\long\def\@savemarbox#1#2{\global\setbox#1\vtop{\hsize\marginparwidth 
  \@parboxrestore\tiny\raggedright #2}}
\newcommand\lref[1]{\ref{#1}%
\@ifundefined{r@DisplaY #1}{}{ (#1)}}

\newcommand\fakelabel[2]{\@bsphack\if@filesw {\let\thepage\relax
   \newcommand\protect{\noexpand\noexpand\noexpand}%
\xdef\@gtempa{\write\@auxout{\string
      \newlabel{#1}{{#2}{\thepage}}}}}\@gtempa
   \if@nobreak \ifvmode\nobreak\fi\fi\fi\@esphack}

\def\SL@margintext#1{{\showlabelsetlabel{\tiny\{\SL@prlabelname{#1}\}}}}
\catcode`\@=12

\def\Empty{}
\newcommand\oplabel[1]{
  \def\OpArg{#1} \ifx \OpArg\Empty {} \else
        \label{#1}
  \fi}
\newtheorem{theoremSt}{Theorem}[section]

\newtheorem{exampleSt}[theoremSt]{Example}
\newtheorem{exerciseSt}[theoremSt]{Exercise}

\newcommand\MakeStEnv[1]{
  \newenvironment{#1}[1]{
  \begin{#1St} \oplabel{##1}%
  \global\def\CrntSt{\thetheoremSt}%
}{ 
  \end{#1St} }
  \newenvironment{#1+}[1]{
  \begin{#1St} \label{##1}%
  \label{DisplaY ##1}%
  \global\def\CrntSt{\thetheoremSt}%
  \def\Labl{##1}\ifx\Labl\Empty{} \else {\em (\Labl)\,}\fi%
}{ 
  \end{#1St} }
}
\MakeStEnv{theorem}
\MakeStEnv{corollary}
\MakeStEnv{proposition}
\MakeStEnv{lemma}
\MakeStEnv{definition}
\MakeStEnv{conjecture}
\MakeStEnv{problem}
\MakeStEnv{question}

\newlength{\saveu}

\newenvironment{pf*}[1]{%
 \begin{proof}[#1]%
}{ 
 \end{proof}
}

\newcommand{\finishproof}[1]{ 
  \def\FPArg{#1}
  \ifx\FPArg\Empty
        \newcommand\FPArg{\CrntSt}  \fi
  \smallbreak\noindent\makebox[\textwidth]{\hfill\fbox{\FPArg}}
  \medbreak\noindent
}

\newcommand\BB{{\mathcal B}}

\newcommand\DD{{\mathcal D}}
\newcommand\EE{{\mathcal E}}
\newcommand\FF{{\mathcal F}}

\newcommand\LL{{\mathcal L}}
\newcommand\MM{{\mathcal M}}
\newcommand\NN{{\mathcal N}}
\newcommand\OO{{\mathcal O}}
\newcommand\PP{{\mathcal P}}

\newcommand\RR{{\mathcal R}}
\newcommand\SSS{{\mathcal S}}

\newcommand\XX{{\mathcal X}}

\newcommand\PMF{{\PP\kern-2pt\MM\FF}}

\newcommand\PML{{\PP\kern-2pt\MM\LL}}

\newcommand\hhat{\widehat}

\newcommand\union{\cup}
\newcommand\intersect{\cap}
\newcommand\bbR{{\mathord{\text{I\kern-2pt R}}}}        
\newcommand\bbH{{\mathord{\text{I\kern-2pt H}}}}        

\newcommand\C{{\mathbb C}}

\newcommand\Z{{\mathbb Z}}
\newcommand\R{{\mathbb R}}
\newcommand\N{{\mathbb N}}

\newcommand\Hyp{{\mathbb H}}

\newcommand\bigrightarrow[1]{\hbox to #1{\rightarrowfill}}
\newcommand\bigleftarrow[1]{\hbox to #1{\leftarrowfill}}

\newcommand\homeo{\cong}
\newcommand\boundary{\partial}
\newcommand\semidir{\mathrel{\hbox{\vrule depth-.03ex height1.1ex\kern-0.15em$\times$}}}

\newcommand\til{\widetilde}

\newcommand\tr{\operatorname{tr}}

\numberwithin{equation}{section}

\catcode`\@=11

\def\subsection{\@startsection{subsection}{2}%
  \z@{.5\linespacing\@plus.7\linespacing}{.5em}%
  {\normalfont\bfseries\centering}}

\def\section{\@startsection{section}{1}%
  \z@{.7\linespacing\@plus\linespacing}{.5\linespacing}%
  {\normalfont\large\bfseries\centering}}

\def\subsubsection{\@startsection{subsubsection}{3}%
  \z@{.5\linespacing\@plus.7\linespacing}{-.5em}%
  {\normalfont\bfseries}}

\catcode`\@=12

\newcommand{\fsubd}{\mathrel{{\scriptstyle\searrow}\kern-1ex^d\kern0.5ex}}
\newcommand{\bsubd}{\mathrel{{\scriptstyle\swarrow}\kern-1.6ex^d\kern0.8ex}}
\newcommand{\fsubeq}{\mathrel{\raise-.7ex\hbox{$\overset{\searrow}{=}$}}}
\newcommand{\bsubeq}{\mathrel{\raise-.7ex\hbox{$\overset{\swarrow}{=}$}}}

\newcommand{\bbar}{\overline}

\newcommand{\tsh}[1]{\left\{\kern-.9ex\left\{#1\right\}\kern-.9ex\right\}}

\newcommand\cv{{\mathcal X}}
\newcommand\PS{{\mathcal{PS}}}
\newcommand\PSLtC{{\operatorname{PSL_2}(\C)}}

\newcommand\PSLtR{{\operatorname{PSL_2}(\R)}}
\newcommand\SLtR{{\operatorname{SL_2}(\R)}}
\newcommand\wh{\operatorname{Wh}}
\newcommand\Isom{\operatorname{Isom}}

\title{On dynamics of $Out(F_n)$ on $\PSLtC$ characters}
\date{\today}

\author{Yair N. Minsky}
\address{Yale University}
\email{yair.minsky@yale.edu}

\thanks{Partially supported by NSF grants DMS-0504019 and DMS-0554321}

\begin{document}

\begin{abstract}
We introduce and study an open set of
$\PSLtC$ characters of a nonabelian free group, on which the action of
the outer automorphism group is properly discontinuous, and which is strictly larger
than the set of discrete, faithful convex-cocompact  (i.e. Schottky)
characters. This implies, in particular, that the
outer automorphism group does not act ergodically on the set of
characters with dense image. Hence there is a difference between the
geometric (discrete vs. dense) decomposition of the characters, 
and a natural dynamical decomposition. 
\end{abstract}

\maketitle

\section{Introduction}
\label{sec:intro}

Let $F_n$ be the free group on $n\ge 2$ generators. Its automorphism group $Aut(F_n)$ 
acts naturally, by precomposition, on $Hom(F_n,G)\equiv G^n$ for any group
$G$. The outer automorphism group $Out(F_n) = Aut(F_n)/Inn(F_n)$  acts
on the quotient $Hom(F_n,G)/G$, where $G$ is understood to act by
inner automorphisms. When $G$ is a Lie group we
consider instead the character variety $\cv(F_n,G) \equiv
Hom(F_n,G)//G$, the quotient in the sense of geometric invariant
theory.

When $G$ is a compact Lie group and $n\ge 3$, Gelander \cite{gelander:outFn}
showed that this action is ergodic, settling a conjecture of Goldman
\cite{goldman:ergodic-su2}, who had proved it for $G= SU(2)$. 
When $G$ is non compact the situation is different because
there is a natural decomposition of $\cv(F_n,G)$, up to measure 0,
into (characters of) {\em dense} and {\em discrete} representations, and in the cases
of interest to us the action on the discrete set is not ergodic,
indeed even has a nontrivial domain of discontinuity.

See Lubotzky \cite{lubotzky:autFn} for a comprehensive survey on the
dynamics of representation spaces, from algebraic, geometric and
computational points of view. 

We will focus on the case of $G=\PSLtC$, where the
interior of the discrete set is the set of Schottky representations. In this
case one can ask if the action is ergodic, or even topologically transitive,
in the set of dense representations, but this turns out to be the
wrong question. In particular: 

\begin{theorem}{larger discontinuity}
There is an open subset of $\cv(F_n,\PSLtC)$, strictly larger than the set of
Schottky characters, which is $Out(F_n)$ invariant, and on which 
$Out(F_n)$ acts properly discontinuously.
\end{theorem}

In other words, the natural ``geometric'' decomposition of
$\cv(F_n,\PSLtC)$, in terms of discreteness vs. density of the image
group, is distinct from the ``dynamical'' decomposition, in terms of
proper discontinuity vs. chaotic action of $Out(F_n)$.

The subset promised in Theorem \ref{larger discontinuity} will be
the set of {\em primitive-stable} representations (see definitions
below). It is quite easy to see that this set is open and $Out(F_n)$
invariant, and that the action on it is properly discontinuous
(Theorems \ref{primitive-stable basics}, \ref{primitive-stable
  disc}). Thus the main content of this note is the observation,
via a lemma of Whitehead on free groups and a little bit of hyperbolic
geometry, that it contains non-Schottky
(and in particular non-discrete) elements. This will be carried out in
Section \ref{sec:whitehead}.

One should compare this with results of Goldman \cite{goldman:action} on
the rank 2 case for $\SLtR$ characters, and with work of Bowditch \cite{bowditch:markoff} on
the complex rank 2 case. Bowditch, and Tan-Wong-Zhang \cite{tan-wong-zhang:markoff}, studied a
condition very similar to primitive stability; we will compare the
two in Section \ref{questions}.

\subsubsection*{Acknowledgements}
Bill Goldman and  Alex Lubotzky got the author interested in this question, and
Lubotzky pointed out the notion of ``redundant representations'', whose negation
leads (eventually) to the primitive-stable condition. Mark Sapir and
Vladimir Shpilrain pointed the author to Whitehead's lemma and
Corollary \ref{primitive blocking}.
The author is also grateful to Dick Canary, Hossein Namazi and Juan Souto
for interesting conversations and perspectives. Juan Souto in
particular pointed out the connection to Otal's work. The referee made
a number of helpful and incisive suggestions.

\section{Background and notation}
\label{sec:background}

In the remainder of the paper we fix $n\ge 2$ and let $G =
\PSLtC$. We also fix a free generating set $X=\{x_1,\ldots,x_n\}$ of $F_n$. 

Note that $Hom(F_n,G)$ can be identified with $G^n$ via $\rho \mapsto
(\rho(x_1),\ldots,\rho(x_n))$ once the generators are fixed. 
The quotient $\cv(F_n) = Hom(F_n,G)//G$, which we denote $\cv(F_n,G)$ or
just $\cv(F_n)$,  is obtained by considering
characters of representations, which in our case corresponds to trace
functions. This quotient naturally has the structure of an algebraic
variety, and differs from the purely topological quotient $Hom(F_n,G)/G$
only at reducible points, i.e. representations whose images fix a
point on $\hat\C$ (see Kapovich \cite{kapovich:book} or Morgan-Shalen
\cite{morgan-shalen:I}). 
Characters of reducible representations form a subset of measure 0
outside all of the sets we shall be considering, and so we shall
be able to ignore the distinction between these two quotients.

\subsubsection*{Geometric decomposition}
There is a natural $Out(F_n)$-invariant decomposition of $\cv(F_n)$ in terms of the geometry
of the action of $\rho(F_n)$ on $\Hyp^3$. Namely, let $\DD(F_n)=\DD(F_n,G)$ denote
the (characters of) discrete and faithful representations, and let
$\EE(F_n)=\EE(F_n,G)$ denote those of representations with dense image in $G$. 

It is fairly well-known (see \cite{breuillard-gelander:dense,breuillard-gelander-souto-storm}) that 
\begin{lemma}{D and E}
$\EE(F_n)$ is nonempty and open, $\DD(F_n)$ is closed, and $\cv(F_n) \setminus
  (\DD(F_n)\union \EE(F_n))$ has measure 0. 
\end{lemma}

The idea of the measure 0 statement is this: If $\rho$ is not faithful
it satisfies some relation; this is a nontrivial algebraic condition,
so defines a subvariety of measure 0. There are a countable number of
such relations. If $\rho$ is not discrete, consider the identity
component of the closure of $\rho(F_n)$ in $G$. This is a connected Lie subgroup
of $G$ and hence is either all of
$G$ (so $\rho(F_n)$ is dense), or solvable (i.e. elementary, fixing a
point in the Riemann sphere), or conjugate to $\PSLtR$. The latter
cases are again detected by algebraic conditions. 

Openness of $\EE(F_n)$ follows from the Kazhdan-Margulis-Zassenhaus
lemma \cite{kazhdan-margulis},  which furnishes a neighborhood $U$ of
the identity in $G$ in which any set of $n$ elements generates either
an elementary or an indiscrete group. 
Generating an elementary group is a nontrivial algebraic condition on
$U^n\subset Hom(F_n,G)$
(for $n\ge 2$), and as above for an indiscrete group not to be dense
is also a nontrivial algebraic condition. Hence an open dense subset
$W$ of $U^n$ consists of dense representations. 
Now given any $\rho$ with dense image, there exist elements
$h_1,\ldots,h_n\in F_n$ such that $(\rho(h_1),\ldots,\rho(h_n))\in W$, 
and since $W$ is open the same is true for $\rho'$ sufficiently close
to $\rho$ -- hence $\rho'$ also generates a dense subgroup.

The fact that $\DD(F_n)$ is closed follows from J\o rgensen's inequality
\cite{jorgensen:inequality}, or 
alternately from the Kazhdan-Margulis-Zassenhaus
lemma. Lemma \ref{D and E} in fact holds for much more
general target groups $G$ -- see 
\cite{breuillard-gelander:dense,breuillard-gelander-souto-storm} for details.

Note, when $G$ is compact $\DD(F_n,G)$ is empty, and in this case Gelander
proved that $Out(F_n)$ acts ergodically on $\cv(F_n,G)$ and hence on
$\EE(F_n,G)$. Our main theorem will show that, in general, the action on
$\EE(F_n,G)$ is not ergodic. 

\subsubsection*{Schottky groups}
A {\em Schottky group} (or representation) $\rho$ is one which is
obtained by a ``ping-pong'' configuration in the sphere at infinity
$\boundary \Hyp^3$. That is, suppose that $D_1,D'_1,\ldots,D_n,D'_n$
are $2n$ disjoint closed (topological) disks in $\boundary\Hyp^3$ and
$g_1,\ldots,g_n\in \PSLtC$ are isometries such that $g_i(D_i)$ is the closure of
the complement of $D'_i$. Then $\{g_1,\ldots,g_n\}$ generate a free
discrete group of rank $n$, called a Schottky group. 
The representation sending $x_i \mapsto
g_i$ is discrete and faithful, and moreover, an open neighborhood of
it in $Hom(F_n,G)$ consists of similar representations. We let $\SSS(F_n)$
denote the open set of all characters of Schottky representations. 

Sullivan \cite{sullivan:rigidity}  proved that 
\begin{theorem}{sullivan}
$\SSS(F_n)$ is the interior of $\DD(F_n)$. 
\end{theorem}
(This theorem is not known to hold for the 
higher-dimensional hyperbolic setting; see \S \ref{questions}.)

To obtain a geometric restatement of the Schottky condition, consider
the {\em limit set} for the action of any discrete group of isometries
on $\Hyp^3$, namely the minimal closed
invariant subset of  $\boundary \Hyp^3$.
The convex hull  in $\Hyp^3$ of this limit set is
invariant, and its quotient by the group is called the {\em convex
  core} of the quotient manifold. The deep usefulness of the convex
core in the study of hyperbolic 3-manifolds was first exploited by Thurston.
The Schottky condition on $\rho$ is
equivalent to the condition that the convex core of $\Hyp^3/\rho(F_n)$
is a compact handlebody of genus $n$ (see Marden
\cite[pp. 98-99]{marden:outercircles} for a discussion and references). 

It is a well-known consequence of the deformation theory of Kleinian
groups (see e.g. Bers \cite{bers:spaces} for background)
that $Out(F_n)$  acts properly discontinuously on
$\SSS(F_n)$: that is, that the set $\{\psi\in Out(F_n): \psi(K)\intersect K \ne
\emptyset\}$ is finite for any compact $K\subset \SSS(F_n)$.
We note that this will also follow from Theorem \ref{primitive-stable disc}.

\section{Primitive-stable representations}
\label{sec:PS}

Let $\Gamma$ be a bouquet of $n$ oriented circles labeled by our fixed
generating set. 
We let $\BB = \BB(\Gamma)$ denote
the set of bi-infinite (oriented) geodesics in  $\Gamma$, as in Bestvina-Feighn-Handel
\cite{bestvina-feighn-handel:I}.  Each such geodesic lifts to
an $F_n$-invariant set of bi-infinite geodesics in $\til
\Gamma$, the universal covering tree (and the Cayley graph
of $F_n$ with respect to $X$).

Let $\boundary F_n$ be the boundary at infinity of
$F_n$, or the space of ends of the tree $\til \Gamma$. We 
have a natural action of $F_n$ on $\boundary F_n$. 
Each element of $\BB$ can be identified with an $F_n$-invariant subset
of $\boundary F_n \times \boundary F_n \setminus \Delta$
(with $\Delta$
the diagonal), i.e. the pairs
of endpoints of its lifts. 
$Out(F_n)$ acts naturally on $\BB$ (and in general on $F_n$-invariant subsets
of $\boundary F_n\times\boundary F_n \setminus \Delta$.). 

Equivalently we can identify $\BB$ with the set of bi-infinite reduced words in the
generators, modulo shift. 
To every conjugacy class $[w]$ in $F_n$ is associated an element of
$\BB$ named $\bbar w$, namely the periodic word determined by
concatenating infinitely many copies of a cyclically reduced
representative of $w$. 

An element of $F_n$ is called {\em primitive} if it is a member of a
free generating set, or equivalently if it is the image of a standard
generator by an element of $Aut(F_n)$. 
Let $\PP=\PP(F_n)$ denote the subset of $\BB$ consisting of $\bbar w$
for conjugacy classes $[w]$ of primitive elements. 
Note that $\PP$ is $Out(F_n)$-invariant.

Given a representation $\rho:F_n\to \PSLtC$ and a basepoint
$x\in\Hyp^3$, there is a unique map $\tau_{\rho,x} : \til \Gamma \to \Hyp^3$
mapping the origin of $\til \Gamma$ to $x$, $\rho$-equivariant, and mapping
each edge to a geodesic segment. Every element of $\BB$ is represented by an
$F_n$-invariant family of leaves in $\til \Gamma$, which map to a family of broken geodesic
paths in $\Hyp^3$. 

\begin{definition}{ps def}
A representation $\rho:F_n\to \PSLtC$ is {\em primitive-stable} if
there are constants $K,\delta$ and a basepoint $x\in\Hyp^3$ such that $\tau_{\rho,x}$
takes all leaves of $\PP$ to $(K,\delta)$-quasi
geodesics. 
\end{definition}

Note that if there is one such basepoint then any basepoint will do,
at the expense of increasing $\delta$. This condition is invariant under
conjugacy and so makes sense for $[\rho]\in \cv(F_n)$. Moreover the
property is $Out(F_n)$-invariant since $\PP$ is $Out(F_n)$-invariant.
Primitive-stability is a strengthening of the
negation of {\em redundancy}, whose relevance was explained to me by
Alex Lubotzky (see \S\ref{questions}).

Let us establish some basic facts.
\begin{lemma}{primitive-stable basics}
\mbox{}
\begin{enumerate}
\item If $\rho$ is Schottky then it is primitive-stable. 
\item Primitive-stability is an open condition in $\cv(F_n)$. 
\item If $\rho$ is primitive-stable then, for every proper free factor $A$
  of $F_n$, $\rho|_A$ is Schottky.
\end{enumerate}
\end{lemma}

\begin{proof}
To see (1), note that if $\rho$ is discrete and faithful then $\tau_{\rho,x}$ is
the lift to universal covers of an embedding (and homotopy-equivalence) of $\Gamma$ into the quotient manifold
$N_\rho = \Hyp^3/\rho(F_n)$. If $\rho$ is Schottky then the
convex core of $N_\rho$ is compact and hence its homotopy-equivalence to the
image of $\Gamma$ lifts to a quasi-isometry of the convex hull of the
group to $\tau_{\rho,x}(\til\Gamma)$.  
It follows that {\em 
  all} leaves in $\BB$ map to uniform quasi geodesics, and in
particular the primitive ones.

Next we prove (2). Let $\rho$ be primitive-stable, 
and fix a basepoint $x$ and quasi-geodesic constants $K,\delta$ as in
the definition.  Let $\tau=\tau_{\rho,x}$. 

Let $L$ be a primitive leaf, with vertex sequence $\{v_i\in\til\Gamma\}$,  and let
$p_i = \tau(v_i)$. The condition that $\tau|_L$ is quasi-geodesic is
equivalent to the following statement: there exist constants $c>0$ and $k\in\N$
such that, if $P_i$ is the hyperplane perpendicularly bisecting the
segment $[p_i,p_{i+k}]$, then for all $j$ $P_{jk}$ separates
$P_{(j-1)k}$ from $P_{(j+1)k}$, and $d(P_{jk},P_{(j+1)k}) > c$. 
This is an easy exercise, and we
note that $K,\delta$ determine $k,c$, and vice versa. 

Now consider a representation $\rho'$ close to $\rho$, and let $\tau'
= \tau_{\rho',x}$. Up to the action of $F_n$ there are only finitely
many sequences of tree edges of length $2k$, and hence the relative
position (i.e. up to isometry) of $P_i$ and $P_{i+k}$, over all primitive leaves and all
$i$, is determined by the $\rho$ image of a finite number of words of $F_n$.
These images each vary continuously with $\rho$, and 
hence for $\rho'$ sufficiently close to $\rho$, we have that the separation and
distance properties for the $P'_i$ still hold, with modified
constants. Hence the primitive leaves are still (uniformly)
quasi-geodesically mapped by $\tau'$.

Finally we prove (3). 
Let $A$ be a proper free factor, so that 
$F_n = A*B$ with $A$ and $B$ nontrivial. 
Suppose $\rho$ is primitive-stable. If $A$ is cyclic, then $\rho|_A$ being Schottky is
equivalent to $A$'s generator having loxodromic image, and this is an
immediate consequence of having a quasi-geodesic orbit. Hence we may
now assume $A$ has rank at least 2. 
By (2), there is a neighborhood $U$ of
$\rho$ consisting of primitive-stable elements. Suppose
$\rho|_A$ were not Schottky. Since $\SSS(A)$ is the interior of 
$\DD(A)$ by Sullivan's theorem \ref{sullivan}, and since $\DD(A)\union
\EE(A)$ is dense in $\cv(A)$ by Lemma  \ref{D and E}, we can perturb $\rho|_A$ arbitrarily
slightly to get a dense representation. Leaving $\rho|_B$ unchanged we 
obtain $\rho'\in U$  with $\rho'|_A$
dense. Now let $g_m\in A$ be an infinite sequence of 
reduced words with $\rho'(g_m)\to id$. 
For any generator $b$ of $B$, a sequence of elementary automorphisms
multiplying $b$ by generators of $A$
(Nielsen moves) takes $b$ to  $g_m b$, which therefore is
primitive. Note that each $g_m b $ is cyclically reduced, so 
primitive-stability of $\rho'$ implies that the axes of $g_m b$
are uniformly quasi-geodesically mapped by $\tau_{\rho',x}$.
But this contradicts the fact that $\rho'(g_m)\to id$ while the length
of $g_m$ goes to infinity.
\end{proof}

Let 
$$
\PS = \PS(F_n)\subset\cv(F_n)$$
be the set of conjugacy classes of primitive-stable
representations. We have shown that $\PS$ is an open $Out(F_n)$
invariant set containing the Schottky set. In fact, 
\begin{theorem}{primitive-stable disc}
The action of $Out(F_n)$ on $\PS(F_n)$ is properly discontinuous. 
\end{theorem}

\begin{proof}
Let $\ell_\rho(g)$ denote the translation length of the geodesic
representative of $\rho(g)$ for $g\in F_n$, and let 
$||g||$ denote the minimal
combinatorial length, with respect to the fixed generators of $F_n$,  of
any element in the conjugacy class of $g$ 
(equivalently it is the word length of $g$ after being cyclically reduced). 

Let $C$ be a compact set in $\PS(F_n)$. For each $[\rho]\in C$ we have a
positive lower bound for $\ell_\rho(w)/||w||$ over primitive elements
of $F_n$, and a continuity argument as in  part (2) of Lemma
\ref{primitive-stable basics} implies that a uniform lower bound
$$
\frac{\ell_\rho(w)}{||w||} > r > 0
$$
holds over all $[\rho]$ in $C$. Now on the other hand an upper bound
on this ratio holds trivially for any $\rho$ by the triangle
inequality applied to any $\tau_{\rho,x}$. Continuity again gives us a
uniform upper bound 
$$
\frac{\ell_\rho(w)}{||w||}  < R
$$
for $[\rho]\in C$
(here, one should choose a compact preimage of $C$
in $Hom(F_n,G)$, which is easy to do). 

Now if $[\Phi]\in Out(F_n)$ satisfies 
$[\Phi](C) \intersect C \ne \emptyset$, we apply the inequalities to
conclude, for $[\rho]$ in this intersection, that 

$$||\Phi(w)|| \le (1/r)\ell_\rho(\Phi(w)) = (1/r)
\ell_{\rho\circ\Phi}(w) \le (R/r)||w||.$$

The proof is then completed by the lemma below. 
\end{proof}

\begin{lemma}{bounded automorphisms}
For any $A$, the set 
$$
\{ f\in Aut(F_n):  ||f(w)|| \le A||w|| \ \ \forall\  \text{primitive}\ w\}
$$
has finite image in $Out(F_n)$. 
\end{lemma}
\begin{proof}
In fact much less is needed; it suffices to have the inequality only for $w$ with
$||w||\le 2$. 
Let $x_1,\ldots,x_n$ be generators of $F_n$, and consider the action of
$F_n$ on the tree $\til \Gamma$ (its Cayley graph). 
For any $i,j\le n$, let $D$ denote the distance between the axis of
$f(x_i)$ and the axis
of $f(x_j)$. A look at the action on the tree indicates, if $D>0$,
that $||f(x_ix_j)|| = 2D + ||f(x_i)|| + ||f(x_j)||$. Hence, since
$||f(x_ix_j)|| \le 2A$, we get an upper bound on $D$. 

An upper bound on the pairwise distances between the axes of the $f(x_i)$ 
implies that there is a point on $\til \Gamma$ which is a bounded
distance from all the axes simultaneously (minimize the sum of
distances to the axes, which is proper and convex unless all the axes coincide). After conjugating $f$ we
may assume that this point is the origin. Now the bound
$||f(x_i)||\le A$ implies a finite number of choices for $f$. 
\end{proof}

\section{Whitehead's lemma and indiscrete primitive-stable representations}
\label{sec:whitehead}

In this section we will prove that $\PS$ is strictly bigger than the
set of characters of Schottky representations. In particular we will define the
notion of a {\em blocking} curve on the boundary of a handlebody, and show

\begin{theorem}{main example}
Let $\rho:F_n\to G$ be discrete, faithful and geometrically finite with one cusp $c$, which is
a blocking curve. Then $\rho$ is primitive-stable.
\end{theorem}

We'll see below (Lemmas \ref{boundary blocking} and \ref{tight masur is
  blocking}) that blocking curves
are a non-empty class. 
Using the deformation theory of hyperbolic 3-manifolds one can then
obtain primitive-stable points in the boundary of Schottky space -- see the end of the section for details.

Note that, since $\PS(F_n)$ is
open, this implies the existence of a rich class of primitive-stable
representations, including dense ones (by Lemma \ref{D and E}), 
as well as
discrete faithful ones with degenerate
ends (since these are topologically generic in the boundary of
Schottky space -- see
\cite{canary-culler-hersonsky-shalen:cusps} and \cite[Cor. A]{canary-hersonsky:dense}).

First let us recall Whitehead's criterion and define the notion of 
{\em blocking}. Whitehead studied the question of which
elements in the free group are primitive. He found a necessary
combinatorial condition, as part of an algorithm that decides the
question of primitivity. 

\subsection*{The Whitehead graph}

As before, we fix a generating set $X=\{x_1,\ldots,x_n\}$ of
$F_n$. For  a word $g$ in
the generators and their inverses, the {\em Whitehead} graph $\wh(g)=\wh(g,X)$
is the graph with $2n$ vertices labeled $x_1,x_1^{-1},\ldots,
x_n,x_n^{-1}$, and an edge from $x$ to $y^{-1}$  for each string
$xy$ that appears in $g$ or in a cyclic permutation of $g$. 
For more information see Whitehead
\cite{whitehead:graph,whitehead:equivalence},
Stallings \cite{stallings:whitehead} and Otal \cite{otal:thesis}.

Whitehead proved 
\begin{lemma}{whitehead} (Whitehead)
Let $g$ be a cyclically reduced word.  
If $\wh(g)$ is connected and has no cutpoint, then $g$ is not primitive. 
\end{lemma}

Define the ``reduced'' Whitehead graph $\wh'(g)$ to be the same as
$\wh(g)$ except that we don't count cyclic permutations of $g$. In
other words we don't consider the word $xy$ where $x$ is the last
letter of $g$ and $y$ is the first, so $\wh'(g)$ may have one fewer
edge than $\wh(g)$. 

Let us say that a reduced word $g$ is {\em primitive-blocking} if it does not
appear as a subword of any cyclically reduced primitive word. 
An immediate corollary of Lemma \ref{whitehead} is:

\begin{corollary}{primitive blocking} 
If $g$ is a reduced word with $\wh'(g)$ connected and without cutpoints, then 
$g$ is primitive-blocking.
\end{corollary}

Let us also say that $g$ is {\em blocking} if some power
$g^k$ is primitive-blocking. A curve on the surface of the handlebody of genus
$n$ is  blocking if a reduced
representative of its conjugacy class in the fundamental group is
blocking (with respect to our
given generators).

An instructive example of a blocking curve occurs for even
rank, when the handlebody is homeomorphic to the product of an
interval with a surface with one boundary component:

\begin{lemma}{boundary blocking}
Let $\Sigma$ be a surface with one boundary component. The curve
$\boundary \Sigma\times\{1/2\}$ in the handlebody $\Sigma\times[0,1]$
is blocking with respect to standard generators of $\pi_1(\Sigma)$; in fact its square is primitive-blocking. 
\end{lemma}

\begin{proof}
Using standard generators $a_1,b_1,\ldots,a_k,b_k$ for $\pi_1(\Sigma)$, the
boundary is represented by $c=[a_1,b_1]\cdots[a_k,b_k]$. $\wh'(c)$ is a
cycle minus one edge (corresponding to $b_k^{-1}a_1$), and $\wh'(c^2)$
contains the missing edge, and so by Whitehead's lemma is blocking. 
\end{proof}

One can construct other blocking curves on the boundary of any
handlebody by explicit games with train tracks. We omit this approach,
and instead study the relationship of the blocking condition to the
Masur Domain in the measured lamination space. 

\subsection*{Laminations and the  Whitehead condition}

Let $\PML(\boundary H)$ denote Thurston's space of projectivized
measured laminations on the boundary of the handlebody \cite{casson-bleiler,travaux}. Within this we
have the Masur domain $\OO$ consisting of those laminations that have
positive intersection number with every lamination that is a limit of
meridians of $H$ \cite{masur:domain}. This is an open set
of full measure in  $\PML(\boundary H)$ \cite{masur:domain, kerckhoff:masurdomain}.

Using Otal's work \cite{otal:thesis}, we can extend Whitehead's condition to laminations on the boundary as
follows. Any free generating set of $F_n$ is dual to a system of disks on
$H$, which cut it into a 3-ball (Nielsen). Given such a system $\Delta
= \{\delta_1,\ldots,\delta_n\}$ of
disks and a lamination $\mu$, realize both $\mu$ and the disk
boundaries in minimal position -- e.g. fix a hyperbolic metric on
$\boundary H$ and make them geodesics. Otal calls $\mu$ {\em tight}
with respect to $\Delta$ if there are no {\em waves} on  $\Delta$ which
are disjoint from $\mu$. A wave on $\Delta$ is an arc properly embedded in
$\boundary H \setminus \Delta$, which is homotopic, rel endpoints,
through $H$ but not through $\boundary H$ into $\Delta$. 
In particular if $\mu$ is tight then no arcs of $\mu\setminus\Delta$
are waves. Hence, 
if  a closed curve is tight with respect to $\Delta$, then
its itinerary through the disks describes a  {\em cyclically reduced}
word in $F_n$ with respect to the dual generators. 

We define $Wh(\mu,\Delta)$ as follows: Cutting along $\Delta$,
$\boundary H$ becomes a planar surface with $2n$ boundary components,
each labeled by $\delta_i^+$ or $\delta_i^-$. Let $w(\mu,\Delta)$
denote the union of the boundary components and the arcs of
$\mu\setminus \Delta$. Then $Wh(\mu,\Delta)$ is obtained from $w(\mu,\Delta)$
by making the boundary components into vertices and identifying
parallel arcs of $\mu$. In particular $Wh(\mu,\Delta)$ is given a
planar embedding. 
If $\mu$ is a single closed
curve and $\Delta$ is dual to the original generators this is
equivalent to the original definition.

Otal proved the following in \cite{otal:thesis}. We give a
proof, since Otal's thesis is hard to obtain. 
\begin{lemma}{tight masur is blocking} {\rm (Otal)}
If $\mu\in\OO(H)$, then there is a generating set with dual disks
$\Delta$ such that $\mu$ is tight with respect to $\Delta$, and $Wh(\mu,\Delta)$ is connected and has no cutpoints.
\end{lemma}

\begin{proof}
First note that $\inf_\delta\{i(\mu,\delta)\}$, where $\delta$ runs
over meridians of $H$, is positive and realized. For if $\{\delta_i\}$ is
a minimizing sequence such that infinitely many of the $\delta_i$ are
distinct then an accumulation point in in $\PML(S)$ will have
intersection number 0 with $\mu$, contradicting $\mu\in\OO(H).$ The
same holds for disk systems, so we may choose a disk system $\Delta$
that minimizes $i(\mu,\Delta)$. 

Now $\mu$ cannot have a wave with respect to
$\Delta$. If it did, then a surgery along such a wave would produce a new $\Delta'$
whose intersection number with $\mu$ is strictly smaller,
contradicting the choice of $\Delta$.

If $Wh(\mu,\Delta)$ is disconnected then there is a loop $\beta$ in
the planar surface $P = \boundary H \setminus \NN(\Delta)$ which
separates the
boundary components, and does not intersect $\mu$ (here $\NN$ denotes
a regular neighborhood). 
This gives a
meridian that misses $\mu$ in $\boundary H$, again contradicting
$\mu\in\OO(H)$. 

If $Wh(\mu,\Delta)$ is connected but has a cutpoint, this is represented by a boundary
component $\gamma\subset \boundary P$, equal to one copy $\delta_i^\pm$ of a component of
$\Delta$ (we are blurring the distinction between the disks in
$\Delta$ and their boundaries in $\boundary H$). Let $\bar\gamma$
denote the other copy. 
Then $\gamma$ separates the planar complex $w(\mu,\Delta)$, and 
we may let $X_1$ be a component of $w(\mu,\Delta) \setminus
\gamma$ which does not contain $\bar\gamma$. Now $X_1\union\gamma$ cuts $P$
into a union of disks-with-holes, one of which, $E$, must contain
$\bar\gamma$.  Let $\beta$ be the boundary component of
$E$ that separates $\bar\gamma$ from $X_1$. 
Any arcs of $w(\mu,\Delta)$ passing through $\beta$ must pass
through $\gamma$ (since $\gamma$ separates the interior of $E$ from
the rest of the graph), but arcs of $X_1$ incident to $\gamma$ do not
meet $\beta$; hence $i(\beta,\mu)$ is strictly less than $i(\gamma,\mu)$.

Because $\beta$ separates $\gamma$ from $\bar\gamma$, 
cutting along $\beta$ and regluing
$\gamma$ to $\bar\gamma$ yields again a connected planar surface -- hence
$\Delta\union\{\beta\}\setminus\{\gamma\}$ is a new disk system, with
strictly smaller intersection number with $\mu$.  Again this is a
contradiction, so we conclude that $Wh(\mu,\Delta)$ is connected and
without cutpoints. 
\end{proof}

Call a lamination $\lambda$ {\em blocking}, with respect to $\Delta$ (or the
dual generators), if $\lambda$ has no waves with respect to $\Delta$,
and there is some $k$ such that every length $k$ subword of the infinite
word determined by a leaf of $\lambda$ passing across $\Delta$ does
not appear in a cyclically reduced primitive word. Note that,
for simple closed curves,  this
coincides with the previous definition of blocking.
An immediate
corollary of the above lemma is:

\begin{lemma}{masur blocking}
A connected Masur-domain lamination
(e.g. a simple closed curve or a filling lamination) on the boundary of a
handlebody is blocking with respect to some generating set.
\end{lemma}

\begin{proof}
Given $\mu\in\OO(H)$ let $\Delta$ be as in Lemma \ref{tight masur is blocking}.
In a connected measured lamination every leaf is dense. Thus a sufficiently
long leaf of $\mu$ would traverse every edge of $\wh(\mu,\Delta)$, and
so the corresponding word is blocking by Corollary \ref{primitive
  blocking}. Note for a  simple closed curve this argument shows that its
square is primitive-blocking. 
\end{proof}

\subsection*{Blocking cusps are primitive-stable}
We can now provide the proof of Theorem \ref{main example}, namely
that a geometrically finite representation with a single blocking cusp
is primitive-stable. 

\begin{proof}[Proof of Theorem \ref{main example}]
Let $N_\rho = \Hyp^3/\rho(F_n)$ be the quotient manifold, and $C_\rho$
its convex core. The geometrically finite hypothesis means that
$C_\rho$ is a union of a compact handlebody $H$ and a subset of a parabolic cusp
$P$ (namely a vertical slab in a horoball modulo $\Z$) along an annulus $A$ with core curve $c$ in
$\boundary H$, which we are further assuming is a blocking curve. 

We will prove that all primitive elements of $F_n$ are represented by
geodesics in a fixed compact set $K\subset C_\rho$. The idea is that
in order to leave a compact set, a primitive element must wind around
the cusp, and this will be prohibited by the blocking property. 

Let $\gamma$ be a closed geodesic in $N_\rho$. Then $\gamma\subset
C_\rho$. The orthogonal projection $P\to\boundary P$ gives a
retraction $\pi:C_\rho\to H$. Let $\hhat\gamma = \pi(\gamma)$. 
Note that if a component $a$ of $\gamma\intersect P$
is long then its image $\pi(a)$
winds
many times around $A$ (in fact the number of times is exponential in
the length of $a$). 

Claim: $\hhat\gamma$ is uniformly quasi-geodesic in $H$, 
with constants independent of $\gamma$. More precisely, the lift
$\til\gamma$ of $\hhat \gamma$ to the universal cover $\til H$ is
uniformly quasi-geodesic with respect to the path metric. This follows
from a basic fact about any family $Q$ of disjoint horoballs  in $\Hyp^3$: 

\begin{lemma}{electro-ambient}
Let $Q$ be a family of disjoint open horoballs in $\Hyp^m$, and let
$\pi_Q:\Hyp^m\to\Hyp^m\setminus Q$ be given by orthogonal projection
from $Q$ to $\boundary Q$ and identity in $\Hyp^m\setminus Q$. If $L$
is a geodesic in $\Hyp^m$ then $\pi_Q(L)$ is a quasigeodesic in
$\Hyp^m\setminus Q$ with its path metric, with constants independent
of $Q$ or $L$. 
\end{lemma}

\begin{proof}
This is closely related to statements in Farb \cite{farb:thesis}
and Klarreich \cite{klarreich:thesis} and can
also be proved in greater generality, e.g. for uniformly separated
quasiconvex subsets of a $\delta$-hyperbolic space.  We will
sketch a proof for completeness. 

Note first, there is a constant $r_0$ such that, if $P_1$ and $P_2$
are horoballs in $\Hyp^m$ with $d(P_1,P_2)\ge 1$, then any two
geodesic segments connecting $P_1$ to $P_2$ in their common exterior lie within
$r_0$-neighborhoods of each other. 

If $Q'\subset Q$ is obtained by retracting horoballs to concentric
horoballs at depth bounded by $r_0$, then it suffices to prove the
theorem for $Q'$. This is because any arc on $\pi_Q(L)$ can be
approximated in a controlled way by an arc on $\pi_{Q'}(L)$. Moreover,
given $L$ it suffices to prove the theorem for the union $Q_L\subset Q$ of horoballs
that intersect  $L$, since $\pi_Q(L) = \pi_{Q_L}(L)$ and 
$\Hyp^m\setminus Q\subset \Hyp^m\setminus Q_L$. We can therefore
reduce to the case that any two components of $Q$ are separated by
distance at least 1, and $L$ peneterates each component of $Q$ to
depth at least $r_0$. 

Let $\beta$ be a geodesic in $\Hyp^m\setminus Q$ with endpoints on
$\pi_Q(L)$. It is therefore a concatentation of hyperbolic geodesics
in $\Hyp^m\setminus Q$ with endpoints on $\boundary Q$, alternating
with geodesics on $\boundary Q$ in the path (Euclidean) metric. Let
$\gamma$ be one of the hyperbolic geodesic segments, with endpoints on
horoballs $P_1,P_2\in Q$. Then $\gamma$ is within $r_0$ of the
component of $L\setminus Q$ connecting these horoballs, and hence
$\gamma$ can be replaced by an arc traveling along $\boundary P_1$,
$L$ and $\boundary P_2$ of uniformly comparable length. Replacing all
hyperbolic segments in this way, and then straightening the arcs of
intersection of the resulting path with $\boundary Q$, we obtain a
segment of $\pi_Q(L)$ whose length is comparable with that of
$\beta. $
Since the endpoints of $\beta$ were arbitrary points of $\pi_Q(L)$,
this gives uniform quasigeodesic constants for $\pi_Q(L)$. 
\end{proof}

Now, $H$ retracts to a spine of the manifold which can be identified
with the bouquet $\Gamma$, so $\til H$ is
naturally quasi-isometric to the  Cayley graph $\til \Gamma$ of
$F_n$. In a tree uniform quasi-geodesics are uniformly close to their
geodesic representatives.

By Lemma \ref{electro-ambient}
$\til\gamma$ is quasi-geodesic in $\til H$, and it
contains a high power of the core curve of $A$. Hence its retraction
to the tree contains a high power of some representative $c'$ of the
curve $c$. This subword $(c')^m$ must be uniformly close to a segment of
the form $(c^*)^m$, where $c^*$ is a cyclically reduced representative
of $c$. It follows that the geodesic representative of $\til\gamma$ in
the tree is uniformly close to a segment of this form as well. Since a
tree is 0-hyperbolic, it follows that the geodesic representative of
$\til\gamma$ in the tree actually contains a high power of $c^*$. 

But since $c$ is blocking, this means that $\gamma$ cannot have been
primitive.  We conclude that all primitive geodesics are trapped in a
fixed compact core of $H$. The retraction of this core to the spine of
$H$ therefore lifts to a quasi-isometry, and stability of
quasigeodesics again implies that each primitive geodesic in the tree
is uniformly quasi-geodesic.
\end{proof}

\begin{proof}[Proof of Theorem \ref{larger discontinuity}]
In view of Lemma \ref{primitive-stable basics} and
Theorem \ref{primitive-stable disc}, all that remains to prove is the
existence of a non-Schottky primitive stable representation. By
Theorem \ref{main example}, any
discrete, faithful and geometrically finite representation with a
single blocking cusp would do.

Given a curve $c$ on the boundary of a handlebody $H$, a sufficent
condition that it be realizeable as the single cusp of a geometrically
finite representation is that it be homotopically nontrivial in the
handlebody, and that any homotopically nontrivial proper annulus in $H$ with a
boundary component on $c$  must be parallel into $c$. 
This is a consequence of
Thurston's geometrization theorem (see e.g. \cite{kapovich:book}) or 
of Maskit \cite{maskit:parabolics}. 

Otal \cite{otal:thesis} shows that any curve in the Masur Domain has
complement satisfying these topological conditions (even more strongly,
that $(H,\boundary H\smallsetminus c)$ is acylindrical, 
so there are no essential annuli at all in $(H,\boundary H\setminus c)$), and hence can
appear as the single cusp of a geometrically finite representation. By
Lemma \ref{masur blocking} such a curve is also blocking with respect 
to some generating set, and hence gives us a primitive-stable
representation. 

When $n$ is even the example of Lemma \ref{boundary blocking}, where $c=\boundary \Sigma$,
also suffices, even though $c$ is not in the Masur domain.
Note that we can explicitly
construct Fuchsian representations  for which $c$ is the unique
cusp, and Lemma \ref{boundary blocking} provides the blocking property.

\end{proof}

\section{Further remarks and questions}
\label{questions}

Having established that $Out(F_n)$ acts properly discontinuously on
$\PS(F_n)$, and that $\PS(F_n)$ is strictly larger than $\SSS(F_n)$,
one is naturally led to study the dynamical decomposition 
of $\cv(F_n)$. In particular we
ask if $\PS$ is the maximal domain of discontinuity,
and what happens in its complement. We have only partial results
in this direction.

\subsection*{Outside $\PS$}

The polar opposite of the primitive-stable characters are the
{\em redundant characters} $\RR(F_n)$, defined (after Lubotzky) as follows: $[\rho]$ is
redundant if there is a proper free factor $A$ of $F_n$ such that
$\rho(A)$ is dense. Note that $\RR(F_n)$ is $Out(F_n)$-invariant. Clearly
$\RR(F_n)$ and $\PS(F_n)$ are disjoint, by Lemma \ref{primitive-stable basics}. 

The set $\EE(F_n)$ of representations with dense image is open (Lemma
\ref{D and  E}), and it follows 
(applying this to the factors) that $\RR(F_n)$ is open. The action of
$Out(F_n)$ on $\RR(F_n)$ cannot be properly discontinuous, and in fact
there is a larger set for which we can show this. 

Let $\PS'(F_n)$ be the set of (conjugacy classes of) representations
$\rho$ which are Schottky on every proper free factor. Hence
$\PS(F_n)\subset\PS'(F_n)$ by Lemma \ref{primitive-stable basics}, and
$\PS'(F_n)$ is still in the complement of $\RR(F_n)$. Let $\RR'(F_n)=
\cv(F_n)\setminus \PS'(F_n)$. 

\begin{lemma}{Rprime nondisc}
\mbox{}
\begin{enumerate}
\item If $n\ge 3$ then $\RR(F_n)$ is dense in $\RR'(F_n)$.
\item No point of $\RR'(F_n)$ can be in a domain of discontinuity for
$Out(F_n)$. Equivalently, any open invariant set in $\cv(F_n)$ on which
$Out(F_n)$ acts properly discontinuously must be contained in
$\PS'(F_n)$. 
\end{enumerate}
\end{lemma}

The proof is quite analogous to part (2) of Lemma
\ref{primitive-stable basics}. In fact we note that part (2) of Lemma
\ref{primitive-stable basics} is an immediate corollary of Lemma
\ref{Rprime nondisc} and Theorem \ref{primitive-stable disc}.

\begin{proof}
For (1), let $[\rho]\in \RR'(F_n)$ and let $A$ be a proper free factor such
that $\rho|_A$ is not Schottky. We may assume $A$ has rank at least 2
since $n\ge 3$. Hence (as in the proof of Lemma
\ref{primitive-stable basics}) $\rho|_A$ is approximated by
representations with dense 
image. It follows that $\rho$ itself is approximated by
representations dense on $A$, so $\RR$ is dense in $\RR'$. 

For (2), 
we will show that for every neighborhood $U$ of $[\rho]\in\RR'(F_n)$ there
is an infinite set of elements $[\phi]\in Out(F_n)$ such that
$[\phi](U)\intersect U\ne \emptyset$. Since $\cv(F_n)$ is locally
compact, this implies $[\rho]$ cannot be in
any open set on which $Out(F_n)$ acts properly discontinuously. 

Suppose first $n\ge 3$. 
Since $\RR(F_n)$ is dense in $\RR'(F_n)$, it suffices to consider the case
that $\rho\in\RR(F_n)$. 
Again let $A$ be a proper free factor on which 
$\rho$ is dense. We can assume that $F_n=A*B$
where $B$ is generated by one element $b$. 
Now let $g_m\in A$ be a sequence such that $\rho(g_m)\to id$, and
let $\phi_m\in Aut(F_n)$
be the automorphism that is the identity on $A$ and sends $b$ to
$g_m b$. Note that $[\phi_m]$ has infinite order in $Out(F_n)$.
The number of powers of $\phi_m$ that take $\rho$ to any fixed
neighborhood of itself goes to $\infty$ as $m\to\infty$, because
$\rho(g_m)\to id$. This concludes the proof for $n\ge 3$. 

If $n=2$, $\RR(F_n)$ is empty. However, every $[\rho]\in\RR'(F_n)$ has
primitive element mapping to a non-loxodromic, so $\rho$ may be
approximated by a representation $\rho'$ sending a generator to an irrational
elliptic. The same argument as above can then be applied to $\rho'$. 
\end{proof}

We remark
that $\PS'$ is indeed strictly larger than $\PS$, at least for even
rank:   If $n=2g$, 
represent the
handlebody as an $I$-bundle over a genus $g$ surface $\Sigma$ with one puncture, and let
$\rho:\pi_1(\Sigma)\to\PSLtC$ be a degenerate surface group with no accidental parabolics
(i.e. every parabolic in $\rho(\pi_1(\Sigma))$ is conjugate to the
image of an element of $\pi_1(\boundary \Sigma)$,
$\rho$ is discrete and faithful, and at least one end of
the resulting manifold is geometrically infinite). The puncture cannot
be in any proper free factor of $F_n$, because it is represented by a
curve on the boundary of the handlebody whose complement is
incompressible. 
Hence the restriction of $\rho$ to any proper free
factor $A$ has no parabolics. By the Thurston-Canary Covering Theorem
\cite{canary:covering} and the Tameness Theorem
\cite{agol:tame,calegari-gabai:tame}, 
$\rho(A)$ cannot be geometrically infinite, so it must be Schottky, and hence
$[\rho]\in \PS'(F_n)$.  On the other hand, every nonperipheral
nonseparating simple curve on $\Sigma$ is primitive, and the ending
lamination of a degenerate end can be approximated by such
curves. Hence the uniform quasigeodesic condition fails on primitive
elements, and $[\rho]$ cannot be in $\PS(F_n)$. 

Note that this example
can be approximated by elements of $\PS(F_n)$, as well as by $\RR'(F_n)$
(the latter by using Cusps are Dense \cite{mcmullen:cusps} or the
Density Theorem for surface groups
\cite{brock-canary-minsky:ELCII}). This leads us to  
ask:

\begin{question}{qn cusps dense}
Is $\PS(F_n)$ the interior of $\PS'(F_n)$?
\end{question}

In particular, in view of Lemma \ref{Rprime nondisc}, a positive
answer would imply 

\begin{conjecture}{domain of discontinuity}
$\PS(F_n)$ is the domain of discontinuity of $Out(F_n)$ acting on
  $\cv(F_n)$. 
\end{conjecture}

We remark that, a priori, there may not be any such domain, i.e. there
may be no maximal set on which the action is properly discontinuous.

\subsection*{Rank 2}

For the free group of rank 2, we have already seen that some of our
statements are slightly different. In particular $\RR(F_2)$ is empty
since no one-generator subgroup of $\PSLtC$ is dense. Moreover,
$\PS'(F_2)$ is exactly the set for which every generator is
loxodromic, and this is dense in $\cv(F_2)$ since it is the complement
of countably many proper algebraic sets. 

Question \ref{qn cusps dense} in particular asks, therefore, whether
$\PS(F_2)$ is dense. It is not clear (to the author) whether this
is true, but there is some evidence against it (see below). 

Another important feature of rank 2 is that the conjugacy class of the
commutator of the generators, and its inverse, are permuted by
automorphisms. It follows that the trace of the commutator is an
$Out(F_2)$-invariant function on 
$\cv(F_2)$, and one can therefore study level sets of this function.

The domain of discontinuity of $Out(F_2)$  was studied by
Bowditch and Tan-Wong-Zhang \cite{tan-wong-zhang:markoff}.  Bowditch defines the following
condition on $[\rho]\in\XX(F_n)$, which Tan-Wong-Zhang call condition BQ:
\begin{enumerate}
\item $\rho(x)$ is loxodromic for all primitive $x\in F_2$.

\item The number of conjugacy classes of primitive elements $x$ such
  that $|\tr(\rho(x))| \le 2$ is finite. 
\end{enumerate}
They show,  using Bowditch's work,  that $Out(F_2)$ acts properly
discontinuously on the invariant open set $BQ$. It is still unclear whether
$BQ$ is the largest such set. 

Note that condition (1) is equivalent to membership in $\PS'(F_2)$. 
It is evident that $\PS(F_2)\subset BQ$, 
and it seems plausible that they are equal.

Note also that computer experiments indicate that the intersection of
$BQ$ with a level set of the commutator trace function is {\em not}
dense in the level set (see 
Dumas \cite{dumas:bear}).
In particular the slice corresponding to trace
$-2$ consists of {\em type-preserving representations} of the
punctured-torus group, i.e. those with parabolic commutator, and
empirically it seems that $BQ$ in this slice coincides with the
quasifuchsian representations (which are all primitive-stable
too, by Theorem \ref{main example} combined with Lemma \ref{boundary
  blocking} ). Bowditch has conjectured that this is in fact the case,
and this seems to be a difficult problem.  At any rate this appears to be
evidence against the density of $\PS(F_2)$.

\subsection*{Ergodicity}
The question of the ergodic decomposition of $Out(F_n)$ on $\cv(F_n)$ 
is
still open. Note, in rank 2, the decomposition must occur along level
sets of the commutator trace function. In rank 3 and higher
our observations indicate that the simplest possible
situation is that, outside $\PS(F_n)$, the action is ergodic, which we
pose as a variation of Lubotzky's original question: 
\begin{question}{qn ergodicity}
Let $n\ge 3$. Is there a decomposition of $\cv(F_n)$ into a domain where the action is
properly discontinuous, and a set where it is ergodic? More pointedly,
does $Out(F_n)$ act ergodically on the complement of $\PS(F_n)$ in
$\cv(F_n)$? 
\end{question}

In Gelander-Minsky \cite{gelander-minsky:redundant} we show that in fact 
the action on $\RR(F_n)$ is ergodic and topologically minimal. 
So if for example $\PS'\setminus \PS$ and $\RR'\setminus
\RR$ have measure 0, we would have a positive answer for the above
question.

\subsection*{Understanding $\PS(F_n)$}
It would also be nice to have a clearer understanding of the boundary
of $\PS(F_n)$, and of which discrete representations $\PS(F_n)$ contains. 

From Lemma \ref{primitive-stable basics} we know that any
discrete faithful representation with cusp curves that have
compressible complement cannot be in $\PS(F_n)$. We've also mentioned the
degenerate surface groups which are in $\PS'(F_n)$ but not $\PS(F_n)$. 

If, however, $\rho$ is discrete and faithful without parabolics and is
not Schottky, then it has an ending lamination which must lie in the
Masur domain, and hence is blocking by Lemma \ref{masur
  blocking}. Hence it would be plausible to expect: 

\begin{conjecture}{degenerate stable}
Every discrete faithful representation of $F_n$ without parabolics is
primitive-stable.
\end{conjecture}

More generally, a discrete faithful representation has a possibly
disconnected ending lamination, whose closed curve components are
parabolics. All the examples we have considered suggest this
conjecture:

\begin{conjecture}{general PS discrete}
A discrete faithful representation of $F_n$ is primitive-stable if and
only if every component of its ending lamination is blocking.
\end{conjecture}

{\em Note:} since the first version of this article appeared, Jeon-Kim
\cite{jeon-kim:PS} gave a positive solution to Conjectures
\ref{degenerate stable} and \ref{general PS discrete}, using the work
of Mj \cite{mitra:CT-kleinian} on Cannon-Thurston maps.

It might also be interesting to think about which representations
with discrete image (but not necessarily faithful) are
primitive-stable. 
In \cite{minsky-moriah:surplus} we construct many
primitive-stable representations whose images uniformize knot
complements. What properties of a marked 3-manifold correspond to
primitive stability? 

Another interesting question is:
\begin{question}{algorithm}
How do we produce computer pictures of $\PS(F_n)$? 
\end{question}
For rank $n=2$, the character variety has complex dimension 2, 
and one can try to draw slices of dimension
1. Komori-Sugawa-Wada-Yamashita developed a program for drawing Bers
slices, which are parts of the discrete faithful locus
\cite{komori-sugawa-wada-yamashita,komori-sugawa:bers}, and Dumas
refined this using Bowditch's work
\cite{dumas:bear}.
In particular what Dumas' program is really doing is drawing 
slices of Bowditch's domain BQ. If indeed $BQ=\PS(F_2)$, 
then this produces images of $\PS(F_2)$ as well.

\subsection*{Other target groups}

The discussion can be extended to other noncompact Lie groups, with
moderate success. Let us consider first the case of
$\Isom_+(\Hyp^d)\homeo SO(d,1)$ for all $d\ge 2$, where $d=3$ is the case
we have been considering. The definition of $\PS$ is unchanged, and
stability of quasigeodesics works in all dimensions in the same
way. Lemmas \ref{D and E} 
Theorem \ref{primitive-stable disc} still hold. 
However, Sullivan's theorem (Theorem \ref{sullivan}) equating Schottky
representations with those in the interior of $\DD(F_n)$ is no longer
available. Schottky representations in higher dimensions can be replaced by
{\em convex-cocompact} representations: discrete and faithful,  with
convex hull of the limit set having a compact quotient. Now the conclusions
of Lemma \ref{primitive-stable basics} must be changed somewhat: A
convex-cocompact representation is certainly still primitive-stable,
but it is not clear that $int(\DD(F_n)) \subset\PS(F_n). $ For a
primitive-stable $\rho$, the proof of Lemma \ref{primitive-stable
  basics} shows that $\rho$ restricted to each proper free factor $A$ is
in $int(D(A))$, but not that it is convex-cocompact.

For $d\ge 3$,
the natural embedding of $\Isom_+(\Hyp^3)$ in $\Isom_+(\Hyp^d)$
clearly preserves primitive-stability and non-discreteness, 
so it is still true that $\PS$ contains indiscrete representations in higher
dimension (and hence dense ones, by Lemma \ref{D and E}).

The case of $d=2$ is slightly trickier. When $n$ is even, we have
given an example of a blocking curve that is the boundary of a
one-holed surface, and so Theorem \ref{main example} shows that a
Fuchsian structure on this surface, which gives an element of
$\cv(F_n,\Isom_+(\Hyp^2)) = \cv(F_n,\PSLtR)$ is primitive-stable but
not Schottky. However when $n$ is odd we have no such example, and it
is unclear to me if $\PS$ contains indiscrete elements. For $n=2$ and $d=2$
Goldman \cite{goldman:action} has described the domain of
discontinuity for the trace $-2$ slice, and proved ergodicity in its complement.

For other noncompact rank-1 semisimple Lie groups, namely isometry groups of
the non-homogeneous negatively curved symmetric spaces, primitive
stability can again be defined in the same way.  The hyperbolic plane
always embeds geodesically in such a space, with its full isometry
group acting. This is not easily extracted from the literature but is
well-known; see Mostow \cite{mostow:monograph}, Bridson-Haefliger
\cite[Chap. II.10]{bridson-haefliger} and Allcock
\cite{allcock:octaveplane} for the requisite machinery, or note that
these spaces are just the real, complex, and quaternionic hyperbolic
spaces, and the Cayley plane, in all of which one can restrict to a
real subspace.  Thus, in these cases our examples for
$F_{2g}$ can be used.

Other cases, such as higher rank semisimple groups, presumably require
a rethinking of the definitions, but there is again no reason to think
that the geometric decomposition of Lemma \ref{D and E} should be the
right dynamical decomposition for the action of $Out(F_n)$. 

A completely different picture holds in the setting of non locally
connected groups. In the
case of $G=SL_2(K)$, with $K$ a non-Archimedean local field of
characteristic $\ne 2$, as well as $G=Aut(T)$ for a tree $T$, 
Glasner \cite{glasner:ergodic-trees}  has shown that $Aut(F_n)$ acts
ergodically on $Hom(F_n,G)$.

\subsection*{Other domain groups}
If we replace $F_n$ by $\pi = \pi_1(S)$ for a closed surface $S$,
$Out(F_n)$ is replaced by $Out(\pi) = MCG(S)$, and the primitive
elements are replaced by their natural analogue, the  simple curves in the surface.
The Schottky representations are replaced by the quasi-Fuchsian
representations $QF(S)$. 
We can define $\PS(\pi)$ in a similar way, but now there is no good reason to
think that  $\PS(\pi)$ is strictly larger than the quasi-Fuschsian
locus. Indeed, every boundary point of  $QF$ can be shown not to lie
in $\PS(\pi)$ (nor in any domain of discontinuity for $MCG(S)$ -- see
also Souto-Storm \cite{souto-storm:mcg} for a related result), and
this is because all parabolics are simple curves, 
and all ending laminations are limits of simple curves. It is
still open as far as I know whether in fact $\PS(\pi)$ is
equal to $QF$; this is closely related to (but formally weaker than)
Bowditch's conjecture in this setting. 

One can also consider  $\cv(H,\PSLtC)$ where $H$ is any
fundamental group of a hyperbolic 3-manifold, and the dynamics are by
$Out(H)$; see Canary-Storm \cite{canary-storm:moduli3}. 
In general, the more complicated the group, the weaker we
should expect the
correspondence between the geometric and dynamical decompositions. An
extreme example is when $H$ is a non-uniform lattice in
$\PSLtC$. In this case $\cv(H)$ has positive dimension, while
Mostow rigidity tells us that $\DD(H)$ is a single point. On the other
hand Mostow also tells us that $Out(H)$ is finite in this case, so
that it acts properly discontinuously on the whole of
$\cv(H)$.

\bibliographystyle{hamsplain}
\bibliography{math}

\end{document}